\documentclass[12pt]{amsart}
\usepackage{amssymb,amscd}
\usepackage[all]{xy}
\title[Kobayashi geodesics]{Kobayashi geodesics in $\sA_g$}
\author[M. M\"oller]{Martin M\"oller}
\address{Max-Planck-Institut f\"ur Mathematik,
Vivatsgasse 7,
53111  Bonn, Germany}
\email{moeller@mpim-bonn.mpg.de  }
\author[E. Viehweg]{Eckart Viehweg}
\address{Universit\"at Duisburg-Essen, Mathematik, 45117 Essen, Germany}
\email{viehweg@uni-due.de}
\thanks{This work has been supported by the DFG-Leibniz program and by the SFB/TR 45
``Periods, moduli spaces and arithmetic of algebraic varieties''.}
 \fussy
\begin{document}
\theoremstyle{plain}
\newtheorem{thm}{Theorem}[section]
\newtheorem{theorem}[thm]{Theorem}
\newtheorem{addendum}[thm]{Addendum}
\newtheorem{lemma}[thm]{Lemma}
\newtheorem{corollary}[thm]{Corollary}
\newtheorem{proposition}[thm]{Proposition}
\theoremstyle{definition}
\newtheorem{remark}[thm]{Remark}
\newtheorem{remarks}[thm]{Remarks}
\newtheorem{notations}[thm]{Notations}
\newtheorem{definition}[thm]{Definition}
\newtheorem{claim}[thm]{Claim}
\newtheorem{assumption}[thm]{Assumption}
\newtheorem{assumptions}[thm]{Assumptions}
\newtheorem{property}[thm]{Property}
\newtheorem{properties}[thm]{Properties}
\newtheorem{example}[thm]{Example}
\newtheorem{examples}[thm]{Examples}
\newtheorem{conjecture}[thm]{Conjecture}
\newtheorem{questions}[thm]{Questions}
\newtheorem{question}[thm]{Question}
\newtheorem{construction}[thm]{Construction}
\numberwithin{equation}{section}
\newcommand{\sA}{{\mathcal A}}
\newcommand{\sB}{{\mathcal B}}
\newcommand{\sC}{{\mathcal C}}
\newcommand{\sD}{{\mathcal D}}
\newcommand{\sE}{{\mathcal E}}
\newcommand{\sF}{{\mathcal F}}
\newcommand{\sG}{{\mathcal G}}
\newcommand{\sH}{{\mathcal H}}
\newcommand{\sI}{{\mathcal I}}
\newcommand{\sJ}{{\mathcal J}}
\newcommand{\sK}{{\mathcal K}}
\newcommand{\sL}{{\mathcal L}}
\newcommand{\sM}{{\mathcal M}}
\newcommand{\sN}{{\mathcal N}}
\newcommand{\sO}{{\mathcal O}}
\newcommand{\sP}{{\mathcal P}}
\newcommand{\sQ}{{\mathcal Q}}
\newcommand{\sR}{{\mathcal R}}
\newcommand{\sS}{{\mathcal S}}
\newcommand{\sT}{{\mathcal T}}
\newcommand{\sU}{{\mathcal U}}
\newcommand{\sV}{{\mathcal V}}
\newcommand{\sW}{{\mathcal W}}
\newcommand{\sX}{{\mathcal X}}
\newcommand{\sY}{{\mathcal Y}}
\newcommand{\sZ}{{\mathcal Z}}
\newcommand{\A}{{\mathbb A}}
\newcommand{\B}{{\mathbb B}}
\newcommand{\C}{{\mathbb C}}
\newcommand{\D}{{\mathbb D}}
\newcommand{\E}{{\mathbb E}}
\newcommand{\F}{{\mathbb F}}
\newcommand{\G}{{\mathbb G}}
\newcommand{\BH}{{\mathbb H}}
\newcommand{\I}{{\mathbb I}}
\newcommand{\J}{{\mathbb J}}
\newcommand{\BL}{{\mathbb L}}
\newcommand{\M}{{\mathbb M}}
\newcommand{\N}{{\mathbb N}}
\newcommand{\BP}{{\mathbb P}}
\newcommand{\Q}{{\mathbb Q}}
\newcommand{\R}{{\mathbb R}}
\newcommand{\BS}{{\mathbb S}}
\newcommand{\T}{{\mathbb T}}
\newcommand{\U}{{\mathbb U}}
\newcommand{\V}{{\mathbb V}}
\newcommand{\W}{{\mathbb W}}
\newcommand{\X}{{\mathbb X}}
\newcommand{\Y}{{\mathbb Y}}
\newcommand{\Z}{{\mathbb Z}}
\newcommand{\rk}{{\rm rk}}
\newcommand{\ch}{{\rm c}}
\newcommand{\Sp}{{\rm Sp}}
\newcommand{\Sl}{{\rm Sl}}
\newcommand{\bsA}{{\overline{\sA}}}
\newcommand{\BEnd}{{\mathbb E}{\rm nd}}
\newcommand{\Mon}{{\rm Mon}}
\newcommand{\MT}{{\rm MT}}
\newcommand{\Hg}{{\rm Hg}}
\newcommand{\End}{{\rm End}}
\newcommand{\Gal}{{\rm Gal}}
\begin{abstract}
We consider Kobayashi geodesics in the moduli space of abelian varieties $A_g$
that is, algebraic curves that are totally geodesic submanifolds for the Kobayashi metric. We
show that Kobayashi geodesics can be characterized as those curves whose logarithmic
tangent bundle splits as a subbundle of the logarithmic tangent bundle
of $A_g$. 
\par
Both Shimura curves and Teichm\"uller curves are examples of Kobayashi geode\-sics, but there
are other examples. We show moreover that non-compact Kobaya\-shi geodesics always
map to the locus of real multiplication and that the $\Q$-irreducibility
of the induced variation of Hodge structures implies that they are defined over a number field.
\end{abstract}
\maketitle
\section*{Introduction}
Let $Y$ be a non-singular projective curve defined over $\C$, let $Y_0\subset Y$ be an open
dense subscheme and let $\W_\Z$ be a polarized $\Z$-variation of Hodge structures on $Y_0$ of weight one. 
So $\W_\Z$ gives rise to a morphism $\varphi_0:Y_0\to \sA_g$, where $g=\frac{\rk(\W)}{2}$. 
\par
It is our aim to characterize those $\varphi$ that are totally geodesic submanifolds
for the Kobayashi metric, in the sequel referred to as Kobayashi geodesics. Roughly speaking,
the characterization will be in terms of the irreducible direct factors of $\W=\W_\Z\otimes_\Z\C$ or 
in terms of the splitting of the natural map 
$$
\varphi^*\Omega^1_{\bsA_g}(\log S_{\bsA_g})\longrightarrow \Omega^1_Y(\log (Y \setminus Y_0)).
$$
\par
We now explain this characterization and the above notation in more detail and then why
the main theorem is a unified treatment of the characterization of Shimura curves and
Teichm\"uller curves - and beyond.
\par
All the conditions we are interested in are preserved by  \'etale coverings. Hence we will allow 
to replace $Y_0$ by an \'etale covering and $\W_\Z$ by its pullback. 
By abuse of notations $\sA_g$ will denote a fine moduli scheme of polarized abelian varieties 
with a suitable level structure. $\bsA_g$ will always denote a {\em Mumford compactification}, i.e. is one of the toroidal compactifications, constructed by Ash, Mumford, Rapoport and Tai and studied by 
Mumford in \cite{mu77}. We write $S_{\bsA_g}=\bsA_g\setminus \sA_g$ for the boundary
and $\varphi:Y \to \bsA_g$ for the extension of $\varphi_0:Y_0\to \bsA_g$ to $Y$.
\par
One of the characterizations we are heading for uses the logarithmic Higgs bundle of a polarized complex variation of Hodge structures (\cite{Si90}, recalled in Section \ref{Kg}).
Writing  $(E,\theta)=(E^{1,0}\oplus E^{0,1}, \theta)$ for the logarithmic Higgs bundle of an 
irreducible direct factor $\V$ of $\W$, and we will say that $(E,\theta)$ (or $\V$) is
{\em maximal Higgs} if $\theta:E^{1,0}\to E^{0,1}\otimes \Omega^1_Y(\log S)$ is an isomorphism.
\par
Then Theorem~\ref{main} states that $\varphi_0:Y_0 \to \sA_g$ is a Kobayashi geodesic
if and only if {\em at least one} of the direct factors $\V$ of $\W$ is maximal Higgs. 
\par
The decomposition of $\W$ in irreducible direct factors is in fact defined over $\bar{\Q}$, hence induced by a decomposition of $\W_{\bar{\Q}}=\W_\Z\otimes_\Z \bar{\Q}$. If one assumes in addition that the general fibre $f_0:X_0\to Y_0$ is simple, then for suitable elements 
$$
\sigma_1={\rm id},\sigma_2,\ldots,\sigma_\ell
$$ 
of the Galois group of $\bar{\Q}$ over $\Q$ one can write
$$
\W=\V\oplus \V^{\sigma_2} \oplus \cdots \oplus \V^{\sigma_\ell},
$$
where $\V^{\sigma_i}$ denotes the conjugate of $\V$ under $\sigma_i$.
\par
The condition ``maximal Higgs'' is not compatible with Galois conjugation.
In particular the Higgs field of $\V^{\sigma_i}$ might be zero or equivalently $\V^{\sigma_i}$
might be unitary. By \cite{VZ04} $\varphi_0:Y_0\to \bsA_g$ is a Shimura curve, hence
a totally geodesic submanifold for the Bergman-Siegel metric, if and only if
{\em all} direct factors of $\W$ are either maximal Higgs or unitary.
\par
Geodesics for the Kobayashi metric have been considered by the first named author 
in \cite{Moe04} under the additional assumption, that $f_0:X_0 \to Y_0$ is a family of Jacobians of a smooth family of curves. In this case $\varphi_0(Y_0)$ is a geodesic for the Kobayashi metric, if and only if the image of $Y_0$ in
the moduli scheme $M_g$ of curves of genus $g$ with the right level structure is a geodesic for the Teichm\"uller metric, hence if and only if $Y_0$ is a Teichm\"uller curve. In particular
$Y_0$ will be affine.

By \cite{Moe04} $Y_0$ being a Teichm\"uller curve is equivalent
to the existence of one direct factor $\V$ of $\W$ which is maximal Higgs.
This implies that it is of rank two, and that its conjugates are
neither unitary nor maximal Higgs. So $\V$ is the variation $\BL$ of Hodge structures,
defined in Section \ref{Kg} starting from logarithmic theta characteristic.
The variation $\BL$ is defined over a totally real number field, and it looks like the
uniformizing variation of Hodge structures on a modular curve, except that it has no $\Q$ structure. Adding the condition that the general fibre of $f_0$ is simple,
the Teichm\"uller curve is determined by the Weil restriction of $\BL$.
In a quite sloppy way one could say that the variation of Hodge structures on a Teichm\"uller curve is like the one on ``a modular curve without a $\Q$-structure''.

The curves $\varphi_0:Y_0\to \bsA_g$ satisfying the equivalent conditions of Theorem
\ref{main} include Teichm\"uller curves as well as Shimura curves, and again they could be seen as a generalization of Shimura curves, obtained by dropping the condition that the variations of Hodge structures are defined over $\Q$. 

Examples of geodesics for the Kobayashi metric are given by curves $Y_0$ on Hilbert modular varieties $U$. By definition, the universal covering $\widetilde{U}$ is a product of disks $\Delta$, and the 
property 'geodesic' just means that one of the projections $\widetilde{U} \to \Delta$
induces an isometry between the universal covering of $Y_0$ and $\Delta$. In Section
\ref{affine} we will show, that all affine geodesics for the Kobayashi metric are obtained in this way, again a result which is well known for Teichm\"uller curves.

In Section \ref{fam} we consider families of curves $\sY_0\to U$ in $\sA_g$, whose general fibre is a Kobayashi geodesic. Assuming that the corresponding variation of Hodge structures is irreducible over $\Q$, or equivalently that the general fibre of the induced family of abelian varieties is simple, we show that such families have to be locally analytically products $Y_0\times U$, hence just the map to $\sA_g$ might vary.
In the affine case, the latter can not happen, which implies that affine Kobayashi geodesics are defined over number fields. In \cite{Moe04} a similar result was stated for Teichm\"uller curves, but its proof, as pointed out by C. McMullen, was incomplete. \vspace{.1cm}
\par
The final version of this article was written during a visit of the second named author to the I.H.E.S., Bures sur Yvette. He would like to thank the members of the Institute for their hospitality. The first named author thanks
the MPIM, Bonn, for its support.

\section{A characterization of Kobayashi geodesics}\label{Kg}
As in the introduction $\sA_g$ denotes a fine moduli scheme of polarized abelian varieties with a suitable level structures, $f_0:X_0\to Y_0$ is a family induced by a morphism 
$\varphi_0:Y_0\to \sA_g$ and $\W_\Z$ is the induced polarized variation of Hodge structures.
We denote the complement $Y \setminus Y_0$ by $S$.
\par
If $K$ is a subfield of $\C$ we write $\W_K=\W_\Z\otimes_\Z K$ and usually
$\W=\W_\C$. Replacing $Y_0$ by an \'etale covering, we may assume that the local monodromy operators around $s\in S$ are unipotent and in addition that $\deg(\Omega^1_Y(\log S))$ is even, hence that there exists a {\em logarithmic theta-characteristic}, i.e. an invertible sheaf $\sL$ on $Y$ with $\sL^2\cong \Omega^1_Y(\log S)$.

For a $\C$-subvariation of Hodge structures $\V\subset \W$ let $\sV$ denote the Deligne extension of
$\V\otimes \sO_{Y_0}$ to $Y$ and let $(E,\theta)$ be the induced logarithmic Higgs bundle, i.e.
the graded bundle $E=E^{1,0}\oplus E^{0,1}$ with respect to the $\sF$-filtration, together with the
Higgs field $\theta:E\to E\otimes \Omega^1_Y(\log S)$ induced by the Gauss-Manin connection.
So $\theta|_{E^{0,1}}=0$ and $\theta(E^{1,0})\subset E^{0,1}\otimes \Omega^1_Y(\log S)$ and we
may consider $\theta$ as a morphism 
$$
\theta:E^{1,0}\to E^{0,1} \otimes \Omega^1_Y(\log S).
$$  
Writing for a coherent sheaf $\sE$ on $Y$ 
$$
\mu(\sE)=\frac{\deg(\sE)}{\rk(\sE)},
$$
the Arakelov inequality, due to Faltings and Deligne (see \cite{Fa83} and \cite{De87}), says that
$$
\mu(\V):=\mu(E^{1,0}) - \mu(E^{0,1}) \leq  \mu(\Omega^1_Y(\log S))=\deg(\Omega^1_Y(\log S)).
$$
\begin{lemma}[\cite{VZ04}]\label{arakeq} The following conditions are equivalent for a non-unitary irreducible subvariation $\V$ of Hodge structures in $\W$.
\begin{enumerate}
\item[a.] $(E,\theta)$ (or $\V$) satisfies the {\em Arakelov equality}, i.e. $\mu(\V) = \mu(\Omega^1_Y(\log S))$.
\item[b.] $(E,\theta)$ is maximal Higgs.
\item[c.] $E^{1,0}$ and $E^{0,1}$ are both stable, and $\theta$ is a morphism between polystable bundles of the same slope.
\item[d.] There exists a rank two, weight one variation of Hodge structures $\BL$ on $Y_0$ and an irreducible unitary local system $\U$ on $Y_0$, regarded as a variation of Hodge structures of bidegree $(0,0)$ such that:
\begin{enumerate}
\item[i.] $\V=\BL\otimes \U$.
\item[ii.] The Higgs bundle of $\BL$ is of the form 
$$
\begin{CD}
\big(\sL\oplus \sL^{-1}, \tau: \sL@> \cong >> \sL^{-1}\otimes \Omega_Y^1(\log S)\big),    
\end{CD}
$$
where $\sL$ is a logarithmic theta characteristic and $\tau$ the induced isomorphism.
\end{enumerate}
\end{enumerate} 
\end{lemma}

As in the introduction $\bsA_g$ denotes a Mumford compactification of $\sA_g$. In particular
$\bsA_g$ is non-singular and the boundary $S_{\bsA_g}$ a normal crossing divisor. Moreover,
as shown in \cite[\S 3]{mu77}, the sheaf $\Omega^1_{\bsA_g}(\log S_{\bsA_g})$ is nef.
We will exclude isotrivial families, hence assume that $\varphi:Y\to \bsA_g$ is finite. 

\begin{definition}\label{Kobdef} Let $M$ be a complex domain and $W$ be a subdomain. $W$
is a {\em totally geodesic submanifold for the Kobayashi metric} if 
the restriction of the Kobayashi metric on $M$ to $W$ coincides with the
Kobayashi metric on $W$. If $W = \Delta$ we call $\Delta$ simply
{\em a (complex) Kobayashi geodesic}.
\par 
We call a map $\varphi_0: Y_0 \to A_g$ {\em a Kobayashi geodesic}, if
its universal covering map 
$$
\tilde{\varphi}_0: \widetilde{Y_0} \cong \Delta \to \BH_g
$$ 
is a Kobayashi geodesic. In particular a Kobayashi geodesic will always be
one-dimensional.
\end{definition}

\begin{theorem} \label{main} Keeping the notations and assumptions introduced above
the following conditions are equivalent:
\begin{enumerate}
\item[a.] $\varphi_0: Y_0 \to \sA_g$ is Kobayashi geodesic.
\item[b.] The natural map $\varphi^* \Omega^1_{\bsA_g}(\log S_{\bsA_g}) \to \Omega^1_Y(\log S)$
splits.
\item[c.] $\W$ contains a non-unitary irreducible subvariation of Hodge structures $\V$ which satisfies the Arakelov equality.
\end{enumerate}
\end{theorem}
Remark that by definition $\varphi_0: Y_0 \to A_g$ is a Kobayashi geodesic if and only if
$Y_0\to \varphi_0(Y_0)$ is \'etale and $\varphi_0(Y_0)\subset \sA_g$ a Kobayashi geodesic. 
Since the conditions b) and c) in Theorem~\ref{main} also are invariant under \'etale
coverings we may always assume that $\varphi_0$ is an embedding.

The proof that c) implies a) is quite easy, and it will be given in Section~\ref{multi}.
There, assuming that $Y_0 \to \sA_g$ is Kobayashi geodesic, we will also construct  
a candidate for the splitting in b) over $Y_0$. In Section~\ref{splittings} we will show that
the splitting extends to a splitting of the sheaves of log-differential forms, finishing the proof that
a) implies b). 

We start in Section \ref{slopes} with the algebraic part of this note, i.e. by showing that b) implies c). In fact, since it hardly requires any additional work, we
will show directly that b) and c) are equivalent.

\section{Slopes of Higgs bundles}\label{slopes}

Instead of studying the log differentials on $Y_0$ we consider the dual sheaves $T_Y(-\log S)=(\Omega_Y^1(\log S))^\vee$ and  $T_{\bsA_g}(-\log S_{\bsA_g})=(\Omega_{\bsA_g}^1(\log S_{\bsA_g}))^\vee$ in this section, and the dual Higgs field
$$
\theta': E^{1,0}\otimes T_Y(-\log S) \longrightarrow E^{0,1}.
$$
Recall that, writing $(E',\theta')$ for the logarithmic Higgs bundle of $\W=\W_Z\otimes_\Z\C$, the pullback of the logarithmic tangent sheaf of $(\bsA_g, S_{\bsA_g})$ to $Y$ is given by 
$$
S^2(E'^{0,1})\subset \sH om(E'^{1,0},E'^{0,1}),
$$
where one uses the polarization to identify $E'^{0,1}$ and the dual of $E'^{1,0}$.
We write $\W=\V_0\oplus \V_1 \oplus \cdots \oplus \V_\ell$, where $\V_0$ is a unitary subsheaf of $\W$ and where $\V_1,\ldots,\V_\ell$ are irreducible non-unitary direct factors.
The logarithmic Higgs bundle of $\V_j$ will be denoted by $(E_j=E_j^{1,0}\oplus E_j^{0,1}, \theta_j)$. The dual Higgs field on $\bsA_g$ induces an identification
\begin{gather}\notag
\varphi^*T_{\bsA_g}(-\log S_{\bsA_g})=S^2(E^{0,1}_0\oplus \cdots \oplus E^{0,1}_\ell)
\cong \big(\bigoplus_j S^2( E^{0,1}_j)\big)\oplus \big(\bigoplus_{i<j} E^{0,1}_i\otimes E^{0,1}_j\big)\\ \label{eqdec}
\subset \sH om(E^{1,0}_0\oplus \cdots \oplus E^{1,0}_\ell,E^{0,1}_0\oplus \cdots \oplus E^{0,1}_\ell)
\cong \bigoplus_{i,j} {E^{1,0}_i}^\vee \otimes E^{0,1}_j. 
\end{gather}
Since $(E_j=E_j^{1,0}\oplus E_j^{0,1}, \theta_j)$ is a Higgs subbundle the composition
\begin{equation}\label{homs}
T_Y(-\log S) \longrightarrow \varphi^*T_{\bsA_g}(-\log S_{\bsA_g})\longrightarrow {E^{1,0}_i}^\vee \otimes E^{0,1}_j
\end{equation}
is zero, except for $i=j\neq 0$. The dual ${E^{1,0}_j}^\vee$ is the complex conjugate of $E^{1,0}_j$ hence isomorphic to $E^{0,1}_\iota$ if $\V_\iota=\bar{\V}_j$.
In this case the right hand side of \eqref{homs} is identified with
$E^{0,1}_\iota \otimes E^{0,1}_j$. If this holds for $\iota=j$, hence if $\V_j$ is defined over $\R$, the map in \eqref{homs} factors through $S^2(E^{0,1}_j)$.  

In order not to be forced to distinguish different cases, we define
for each of the non-unitary irreducible $\C$-subvariation of Hodge structures $\V_j$ in
$\W=R^1f_*\C_{X_0}$, with Higgs bundle $(E_j,\theta_j)$,  
\begin{equation}\label{decomp}
\sT_j=\begin{cases} 
S^2(E_j^{0,1}) & \mbox{if } \V_j \mbox{ is defined over }\R\vspace{.1cm}\\
{E_j^{1,0}}^\vee\otimes E_j^{0,1} & \mbox{otherwise. }
\end{cases}
\end{equation}
The sheaf $\sT_j$ is a direct factor of $\varphi^*T_{\bsA_g}(-\log S_{\bsA_g})$, uniquely determined by $\V_j$. As we have seen above, the image of
$T_Y(-\log S)$ lies in $\sT_1 \oplus \cdots \oplus \sT_\ell$.

As in Theorem~\ref{main} c) we start with a non-unitary irreducible direct factors of $\W$, say $\V_1$. For simplicity we will drop the lower index ${}_1$ in the sequel. The Higgs bundle $(E,\theta)$ gives rise to a morphism
$$\begin{CD}
\eta=\eta_{(E,\theta)}:T_Y(-\log S) @>>> \sT
\end{CD}$$
and since $\V$ is non-unitary, hence $\theta$ non-trivial, $\eta$ is injective.
\begin{lemma}\label{arakeq-splits}
Let $\V$ be non-unitary and satisfying the Arakelov equality. Then for
the direct factor $\sT$ of $\varphi^* T_{\bsA_g}(-\log S_{\bsA_g})$, induced by $\V$,
the natural injection $T_Y(-\log S) \to \sT$ splits.

In particular the existence of such a direct factor $\V$ of $\W$ implies that
$$
T_Y(-\log S) \longrightarrow \varphi^* T_{\bsA_g}(-\log S_{\sA_g})
$$ 
splits.
\end{lemma}
\begin{proof}
Using the notation from Lemma~\ref{arakeq}, d) one finds that $\V=\BL\otimes \U$, with $\U$ unitary.
So the logarithmic Higgs bundle of $\U$ is of the form $(\sU,0)$. The sheaf
$T_Y(-\log S)\cong \sL^{-2}$ is the subsheaf of 
$$
{E^{1,0}}^\vee\otimes E^{0,1}\cong \sL^{-1} \otimes \sU^\vee \otimes \sL^{-1}\otimes \sU,
$$ 
induced by the homotheties in $\sH om(\sU,\sU)$. In particular one has a splitting of $\eta(T_Y(-\log S))$ as subsheaf of ${E^{1,0}}^\vee\otimes E^{0,1}$, hence of $\sT$. 
\end{proof}
As a next step we will show the converse of the first part of Lemma~\ref{arakeq-splits}. 
\begin{proposition}\label{splits-arakeq}
Assume that $\eta(T_Y(-\log S))$ is non-zero and a direct factor of $\sT$.
Then $\V$ satisfies the Arakelov equality.
\end{proposition}
\begin{proof}
The sheaf $\sT$ is a direct factor of ${E^{1,0}}^\vee\otimes E^{0,1}$, hence by assumption the invertible sheaf $T_Y(-\log S)\cong \eta(T_Y(-\log S))$ as well.

By Lemma~\ref{arakeq} the conclusion of Proposition~\ref{splits-arakeq} 
would imply the stability of the sheaves $E^{1,0}$ and $E^{0,1}$. 
The strategy of the proof of the proposition will be to verify first the semistability of both sheaves. This will give the semistability of ${E^{1,0}}^\vee\otimes E^{0,1}$, which implies in turn that the slope of the direct factor $T_Y(-\log S)$ has to be $\mu(E^{0,1})-\mu(E^{1,0})$.\vspace{.2cm}

Let $0=\sG_0 \subset \sG_1 \subset \cdots \subset \sG_\alpha$ be the Harder-Narasimhan (HN) filtration of ${E^{1,0}}^\vee\otimes E^{0,1}$.
\begin{claim}
For some $\ell>0$ the sheaf $T_Y(-\log S)$ maps to a direct factor of $\sG_{\ell}/\sG_{\ell-1}$.  
\end{claim}
\begin{proof}
Since $\sT$ is a direct factor of ${E^{1,0}}^\vee\otimes E^{0,1}$ we can write
$$
{E^{1,0}}^\vee\otimes E^{0,1}=\sG'\oplus T_Y(-\log S)
$$ 
for some $\sG'$. Let $0=\sG'_0 \subset \sG'_1 \subset \cdots \subset \sG'_{\alpha'}=
\sG'$ be the HN-filtration of $\sG'$. \vspace{.2cm}

If $\mu(T_Y(-\log S)) < \mu(\sG'_\nu/\sG'_{\nu-1})$ define $\sG''_\nu=\sG'_\nu$.\vspace{.2cm} 

If $\mu(T_Y(-\log S)) = \mu(\sG'_\nu/\sG'_{\nu-1})$ choose $\sG''_\beta=\sG'_\beta \oplus T_Y(-\log S)$ for all $\beta\geq \nu$.\vspace{.2cm}

If  $\mu(\sG'_{\nu-1}/\sG'_{\nu-2}) > \mu(T_Y(-\log S)) > \mu(\sG'_\nu/\sG'_{\nu-1})$ choose
$\sG''_\beta=\sG'_{\beta-1} \oplus T_Y(-\log S)$ for all $\beta\geq \nu$.\vspace{.2cm}

\noindent
Obviously $\sG''_\bullet$ is a second HN-filtration of ${E^{1,0}}^\vee\otimes E^{0,1}$, and the uniqueness of the HN-filtration implies the claim.
\end{proof}
Let $\sF^{1,0}_\bullet$ and $\sF^{0,1}_\bullet$ be the HN-filtrations of $E^{1,0}$ and $E^{0,1}$, and write
$$
\mathfrak g \mathfrak r _{\sF^{1,0}_\bullet}=\bigoplus H^{1,0}_\eta, \ \ 
\mathfrak g \mathfrak r _{\sF^{0,1}_\bullet}=\bigoplus H^{0,1}_\eta.
$$
So $\sG_{\nu}/\sG_{\nu-1}$ is the direct sum of copies of ${H^{1,0}_{\eta_1}}^\vee\otimes
H^{0,1}_{\eta_2}$, all of the same slope. Since those sheaves are semi-stable, we find
some pairs $(\eta_1,\eta_2)$ such that $T_Y(-\log S)$ is a direct factor of ${H^{1,0}_{\eta_1}}^\vee\otimes
H^{0,1}_{\eta_2}$. In particular
\begin{equation}\label{eqsplit}
\mu(H^{0,1}_{\eta_2})=\mu(T_Y(-\log S)) + \mu(H^{1,0}_{\eta_1}).
\end{equation}
If one replaces $\V$ by its complex conjugate, the role of ${H^{1,0}_{\eta_1}}^\vee$ and 
$H^{0,1}_{\eta_2}$ is interchanged. So we may assume that
\begin{equation}\label{eqsplit2}
\mu(H^{0,1}_{\eta_2}) \geq \mu({H^{1,0}_{\eta_1}}^\vee)=- \mu(H^{1,0}_{\eta_1}),
\end{equation}
for at least one pair $(\eta_1,\eta_2)$ for which \eqref{eqsplit} holds.

Consider next the restriction of the Higgs field to $\sF^{1,0}_{\eta_1}$ with image
$$
\sB\otimes \Omega^1_Y(\log S)\subset E^{0,1}\otimes \Omega^1_Y(\log S).
$$
Writing $\sK$ for the kernel of $\sF^{1,0}_{\eta_1}\to \sB\otimes \Omega^1_Y(\log S)$, one gets 
$$
\deg(\sB)-\rk(\sB)\cdot \mu(T_Y(-\log S))= \deg(\sB\otimes \Omega^1_Y(\log S))\\
=\deg(\sF^{1,0}_{\eta_1}) - \deg(\sK)
$$
Both, $(\sK,0)$ and $(\sF^{1,0}_{\eta_1} \oplus \sB, \theta|_{\sF^{1,0}_{\eta_1}})$ are Higgs subbundles of $(E,\theta)$. So Simpson's correspondence \cite{Si90} implies that $0\geq \deg(\sK)$, or equivalently
\begin{equation}\label{eqsplit3}
\deg(\sB)+\deg(\sF^{1,0}_{\eta_1}) \geq 2\cdot \deg(\sF^{1,0}_{\eta_1}) + \rk(\sB)\cdot \mu(T_Y(-\log S)),
\end{equation}
and that $0 \geq \deg(\sB)+\deg(\sF^{1,0}_{\eta_1})$. Since obviously
$\rk(\sB) \leq \rk(\sF^{1,0}_{\eta_1})$ and since $\mu(T_Y(-\log S))<0$, this together with the inequality \eqref{eqsplit3} implies that
\begin{equation}\label{eqsplit5}
0 \geq \deg(\sB)+\deg(\sF^{1,0}_{\eta_1}) \geq \rk(\sF^{1,0}_{\eta_1})\cdot (2 \mu(\sF^{1,0}_{\eta_1}) +  \mu(T_Y(-\log S))).
\end{equation}     
By definition of the HN-filtration one has
\begin{equation}\label{eqsplit4}
\deg(\sF^{1,0}_{\eta_1}) \geq \rk(\sF^{1,0}_{\eta_1})\cdot \mu(H^{1,0}_{\eta_1})\mbox{ \ \ or \ \ }
\mu(\sF^{1,0}_{\eta_1}) \geq \mu(H^{1,0}_{\eta_1}),
\end{equation}
and since we assumed \eqref{eqsplit2} one finds that
$2 \mu(\sF^{1,0}_{\eta_1}) \geq\mu(H^{1,0}_{\eta_1}) - \mu(H^{0,1}_{\eta_2})$ or
\begin{multline}\label{eqsplit6}
\rk(\sF^{1,0}_{\eta_1})\cdot (2 \mu(\sF^{1,0}_{\eta_1})  +  \mu(T_Y(-\log S))) \geq\\
\rk(\sF^{1,0}_{\eta_1})\cdot (\mu (H^{1,0}_{\eta_1}) - \mu(H^{0,1}_{\eta_2}) +  \mu(T_Y(-\log S))).
\end{multline}
By \eqref{eqsplit} the left hand side in \eqref{eqsplit6} is zero,
hence \eqref{eqsplit5} and \eqref{eqsplit6} must both be equalities. The first one implies that $\rk(\sB)=\rk(\sF^{1,0}_{\eta_1})$ and that $0 = \deg(\sB)+\deg(\sF^{1,0}_{\eta_1})$. So $\sK=0$ and since we assumed $\V$ to be irreducible, $\sF^{1,0}_{\eta_1} =E^{1,0}$ and $\sB=E^{0,1}$. In addition \eqref{eqsplit4} has to be an equality, so $E^{0,1}$ and $E^{0,1}\cong
E^{1,0}\otimes T_Y(-\log S)$ are both semistable. 

As said already at the beginning of the proof, the semistability implies the Arakelov
equality. In fact we could also refer to the equality in \eqref{eqsplit} since we know that
$$\mu(E^{1,0})-\mu(E^{0,1})=\mu(H^{1,0}_{\eta_1})-\mu(H^{0,1}_{\eta_2}).$$
\end{proof}
\begin{proof}[Proof of ``b) $\Longleftrightarrow$ c)'' in Theorem~\ref{main}]
By Lemma \ref{arakeq-splits} it only remains to show that b) implies c).
Recall that one has a factorization 
$$
\begin{CD}
T_Y(-\log S) @>>>  \sT_1 \oplus \cdots \oplus \sT_\ell
@>>> \varphi^*T_{\bsA_g}(-\log S_{\bsA_g})
\end{CD}
$$
and that the composition with the projection ${\rm pr}:\varphi^*T_{\bsA_g}(-\log S_{\bsA_g}) \to  {E_j^{1,0}}^\vee\otimes E_j^{0,1}$ is induced by the Higgs field of
$\V_j$. 

The splitting $\varphi^*T_{\bsA_g}(-\log S_{\bsA_g}) \to T_Y(-\log S)$ of $T_Y(-\log S)$
defines a splitting $\hat{\eta}:\sT_1 \oplus \cdots \oplus \sT_\ell \to  T_Y(-\log S)$, and the
the composition
$$\begin{CD}
T_Y(-\log S) @> \subset >> \sT_j @> \subset >> \sT_1 \oplus \cdots \oplus \sT_\ell  
@> \hat{\eta} >>  T_Y(-\log S)
\end{CD}$$
is either zero, or an isomorphism. So one finds some $j>0$ for which $T_Y(-\log S)$ is a direct factor of
$\sT_j$ and, by Proposition~\ref{splits-arakeq}, $\V_j$ satisfies the Arakelov equality.
\end{proof}

\section{Geodesics and multidiscs}\label{multi}

We recall some facts needed for a characterization of Kobayashi geodesics in 
$\sA_g$. The universal covering of $\sA_g$ is isomorphic to the Siegel 
half space. For $Z \in  M^{g\times g}(\C)$, we denote bz $||\cdot||$ the
operator matrix norm. Via the Cayley transformation it also has 
a realization  as a bounded symmetric domain
$$ 
\sD_g = \{Z \in M^{g\times g}(\C) \, | \, I_g - ZZ^* > 0  \} 
= \{Z \in M^{g\times g}(\C) \, | \, ||Z|| <1  \} 
\subset \Delta^{g^2} \subset \C^{g^2}
$$
that will be more convenient to work with in the next section
when we discuss boundary components.
The intersection of $\sD_g$ with the diagonal in $\C^{g^2}$
is isomorphic to $\Delta^g$, a totally geodesic submanifold for the 
Bergman metric 
on $\sD_g$. Submanifolds of $\sD_g$ that are isomorphic to $\Delta^r$ for 
some $r$ and totally geodesic for the Bergman metric are 
called {\em multidiscs}.  It is well known that $r\leq g$, and the 
multidiscs of maximal dimension, hence isomorphic to
$\Delta^g$, will be called {\em polydiscs}. Note that the Cayley 
transformation maps diagonal matrices to diagonal matrices and that hence
polydiscs are given by diagonal matrices in $\BH_g$ as well. 
\par
Fix a base point $0\in \sD_g$. The fixgroup of $0$ in the isometry group of $\sD_g$ is a maximal compact subgroup $K$. It acts on the set of polydiscs and 
\begin{equation} \label{eqdisccover}
\sD_g=\bigcup_{k\in K}k(\Delta^g).
\end{equation}
In particular for any $v \in T_{0,\sD_g}$ there is a polydisc in 
$\sD_g$ tangent to $v$, and for a second point $p\in \sD_g$ one can find some $k$ with
$0,p\in k(\Delta^g)$. 
\par
The {\em Carath{\'e}odory metric} $c_M$ on $M$ is defined by
$$ 
c_M(x,y) = \sup_{p: M \to \Delta} \rho(p(x),p(y)),
$$
where $\rho$ is the Poincar{\'e} metric on $\Delta$. One shows 
using the Arzela-Ascoli theorem that this supremum is 
attained for some map $p$.
\par
In the following proposition we summarize a structure theorem for
Kobayashi geodesics in $\sD_g$. Part a) and b) are well-known and part c)
is the special case of a structure result for Kobayashi geodesic
in Hermitian symmetric spaces by Abate (\cite{Ab93}). For the
convenience of the reader, we repeat his proof.
\par
\begin{proposition}\label{interprkob}
Suppose we are given a holomorphic map $\tilde{\varphi}_0: \Delta \to \sD_g$. 
\begin{enumerate}
\item[a.)] If there is a map $p: \sD_g \to \Delta$ such that
$p \circ \tilde{\varphi}_0$ is an isometry for the Kobayashi metric
on $\Delta$, then  $\tilde{\varphi}_0$ is a Kobayashi geodesic.
\item[b.)] Conversely, if $\tilde{\varphi}_0$ is a Kobayashi geodesic, 
then there is a holomorphic map $p: \sD_g \to \Delta$ such that 
$p \circ \tilde{\varphi}_0 = {\rm id}$.
\item[c.)] More precisely, if $\tilde{\varphi}_0$ is a Kobayashi geodesic and $0 \in \Delta$
a base point with $\tilde{\varphi}_0(0)=0\in \sD_g$, then 
there exists an integer $1 \leq \ell \leq g$, a matrix $U \in U(g)$,
a morphism $\tilde{\psi}_0: \Delta \to \sD_{g-\ell}$, and a totally geodesic embedding 
$$
i:\Delta^\ell \times \sD_{g-\ell}\longrightarrow \sD_g
$$ 
for the Bergman metric with
\begin{enumerate}
\item[i.)] $\displaystyle i_0=i|_{\Delta^\ell\times \{0\}}:\Delta^\ell \longrightarrow \sD_g$ is 
a multidisk.
\item[ii.)] $U\cdot \tilde{\varphi}_0 \cdot {}^tU = i\circ ({\rm diag}, \tilde{\psi}_0)$,
where $\rm diag$ denotes the diagonal embedding. 
\item[iii.)] $\tilde{\psi}_0(\tau)$ is a matrix $Z_{g-\ell}(\tau)$ with $||Z_{g-\ell}(\tau)|| < |\tau| $.
\end{enumerate}
\end{enumerate}
\end{proposition}
\par
\begin{proof}
Part a) is immediate from the distance-decreasing
property of the Kobayashi metric.
\par
For c) fix some point $\tau_0 \in \Delta \setminus \{0\}$. By \eqref{eqdisccover}
there exists some $U \in U(g)$ such that $M=U\tilde{\varphi}_0(\tau_0){}^tU$
is a diagonal matrix. Let $m_{ii}$ denote the diagonal entries
of $M$ and let $\ell$ be the number of entries with 
$d_\Delta(m_{ii},0) = d_\Delta(\tau_0, 0)$. Since $\varphi_0$
is a Kobayashi geodesic, $\ell$ is not zero. 

Modifying $U$ we may assume that $d_\Delta(m_{ii},0) = d_\Delta(\tau_0, 0)$ for $i=1,\ldots,\ell$, and
multiplying $U$ with a product of matrices in $U(1)\times \cdots \times U(1)$, we may moreover suppose 
that $m_{ii} = \tau_0$ for $i=1,\ldots,\ell$. Since we will need this later, let us state 
to which extent $U$ is unique: 
\begin{enumerate}
\item[]
\begin{enumerate}
\item[iv.)]
If for some point $\tau_0 \in \Delta \setminus \{0\}$, some $\ell' \leq \ell$ and some
$U$ the map $U \cdot \tilde{\varphi}_0 \cdot {}^tU$ has the block form  
$$
(U \cdot \tilde{\varphi}_0(\tau_0) \cdot {}^tU) = \left(\begin{matrix}
A & 0 \\
0 & B 
\end{matrix}
\right),
$$ 
with $A=(a_{ij})$ an $\ell' \times \ell' $ diagonal matrix with $d_\Delta(a_{ii},0) = d_\Delta(\tau_0, 0)$, then 
$$
U = \left(\begin{matrix}
U_{\ell'} & 0 \\
0 & U_{g-\ell'}
\end{matrix}
\right),
$$ 
with $U_{g-\ell'}\in U(g-\ell')$ and with $U_\ell' \in U(1)\times \cdots \times U(1)$.
\end{enumerate}
\end{enumerate}
For some matrices $Z_{11} \in M^{\ell \times \ell}(\C)$, $Z_{12} \in M^{\ell \times (g-\ell)}(\C)$,
$Z_{22} \in M^{(g-\ell) \times (g-\ell)}(\C)$, and $Z_{21} \in M^{(g-\ell) \times \ell)}(\C)$
write 
$$ U\tilde{\varphi}_0(\tau){}^tU = \left(\begin{matrix}
Z_{11}(\tau) & Z_{12}(\tau) \\
Z_{21}(\tau) & Z_{22}(\tau) 
\end{matrix}
\right). $$
We now check that 
$$\psi: \tau \mapsto (Z_{11}(\tau), Z_{12}(\tau))$$
maps to $B_\ell$, the ball of radius $\ell$ in $M^{\ell\times g}(\C)$.
In fact, $1_g - \tilde{\varphi}_0(\tau)\tilde{\varphi}_0(\tau)^*$ is positive definite, 
thus $1_\ell - Z_{11}(\tau)Z_{11}^*(\tau) - Z_{12}(\tau)Z_{12}^*(\tau)$
is positive definite. This implies that the trace of this matrix is
positive, or equivalently, writing $(Z_{11}(\tau), Z_{12}(\tau))=(z_{ij})$ that
$$ \sum_{i=1}^g \sum_{j=1}^\ell |z_{ij}(\tau)| < \ell, $$
which implies the claim.
\par
By construction, $d_{B_\ell}(\psi(\tau_0),0) = d_\Delta(\tau_0,0)$.
Since $B_\ell$ is a strictly convex bounded domain, $\psi:\Delta \to B_\ell$ is a 
Kobayashi geodesic by \cite{Ve81} Proposition~3.3. Moroever, by
\cite{Le81} Theorem~2  there is a unique Kobayashi geodesic 
in $B_\ell$ with $\psi(\tau_0) = 
(\tau_0 \cdot I_d, 0)$, and by \cite{Ve81}, Example on p.~386, it is linear, hence
given by $\psi (\tau) = (\tau \cdot 1_\ell,0)$. In particular $Z_{12}(\tau)=0$ and
by the symmetry of the period matrix one finds that $Z_{21}(\tau)=0$ as well. 

The maps $p$ claimed to exist in part b) are just the projection maps to
$Z \mapsto z_{ii}$, for one $i \in \{1,\ldots, \ell\}$.
\end{proof}
Let us return to the situation considered in Section~\ref{Kg}. So the
map $\varphi_0:Y_0\to \sA_g$ is induced by the variation of Hodge structures $\W_\Z$,
and the image $\tilde{\varphi}_0(\widetilde{Y}_0)$ of the universal covering map is isomorphic to a disc. 

The variation of Hodge structures $\W$ is given by a representation of $\pi(Y_0,*)$, 
hence by a homomorphism $\rho:\pi(Y_0,*) \to \Sp(2g,\R)$. Its image lies in the subgroup 
$$
H = \{ \sigma \in \Sp(2g,\R); \ \sigma(\widetilde{Y}_0)= \widetilde{Y}_0 \}.
$$
\begin{corollary}\label{invariant} The multidisc $i_0:\Delta^\ell \to D_g$ is $H$ invariant, and  
the block decomposition of $\tilde{\varphi}_0$ is preserved by the action of $H$. Morover
there is a subgroup $H_0\subset H$ of finite index which respects the factors of $i_0:\Delta^\ell \to D_g$.
\end{corollary}
\begin{proof} We may assume that in Proposition \ref{interprkob} the matrix $U$ is the
identity matrix. The block decomposition $i:\Delta^\ell\times \sD_{g-\ell}\to \sD_g$
in Proposition \ref{interprkob}, c.) induces an embedding of isometry groups
$$
\hat{i}:\Sl(2,\R)\times \cdots \times \Sl(2,\R) \times \Sp(2(g-\ell),\R) \longrightarrow \Sp(2g,\R).
$$
Let $G$ be the subgroup of $\Sp(2g,\R)$ generated by the image of $\hat{i}$ and by the permutation of the
factors of $i_0(\Delta^\ell)$. So we have to show that $H$ is a subgroup of $G$.

Consider for $h\in H$ the map $h \circ \tilde{\varphi}_0$. Since $h$ leaves $\widetilde{Y}_0$ invariant,
$$
h\circ \tilde{\varphi}_0 (0)=p \in  i(\Delta^\ell\times \sD_{g-\ell}),
$$
and there exists some element $g$ in the image of $\hat{i}$ such that $g(p)=0$. So it is sufficient to show that
$g\circ h$ lies in the image of $G$. To this aim, we apply Proposition \ref{interprkob}, c.) to the
Kobayashi geodesic curve
$$
\tilde{\varphi}'_0=g \circ h \circ \tilde{\varphi}_0:\Delta \longrightarrow \sD_g.
$$
Since $h$ leaves $\widetilde{Y}_0$ invariant, for a point $\tau_0\in \Delta\setminus \{0\}$ 
the image $h\circ \tilde{\varphi}_0(\tau_0)$ is of the form
$$
\left(\begin{matrix}
A & 0 \\
0 & B 
\end{matrix}
\right)
$$
where $A=(a_{ij})$ is an $\ell\times \ell$ diagonal matrix. By the choice of $g$ the same holds for
$\tilde{\varphi}'_0(\tau_0)$. Moreover, since $h$ and $g$ are isometries, 
$d_\Delta(a_{ii},0)=d_\Delta(\tau_0,0)$ we may choose by the condition iv.), stated in the proof of
Proposition \ref{interprkob}, a matrix 
$$
U = \left(\begin{matrix}
U_{\ell} & 0 \\
0 & U_{g-\ell}
\end{matrix}
\right),
$$ 
with $U_{g-\ell}\in U(g-\ell)$ and with $U_\ell \in U(1)\times \cdots \times U(1)$, for which the conclusion of
Proposition \ref{interprkob}, c) holds with $\tilde{\varphi}_0$ replaced by $\tilde{\varphi}'_0$.
So for some $\ell' \leq \ell$ one has 
$$
U\cdot \tilde{\varphi}_0 \cdot {}^tU = i'\circ ({\rm diag}, \tilde{\psi}'_0)\mbox{ \ \  for \ \ }
i':\Delta^{\ell'} \times \sD_{g-\ell'} \longrightarrow \sD_g.
$$
Thus, writing by abuse of notations the conjugation on the other side 
and $Z_{g-\ell'}(\tau)$ for $\tilde{\psi}'_0(\tau)$, we obtain
$$
\tilde{\varphi}'_0(\tau) =\left(\begin{matrix}
\tau \cdot 1_\ell' & 0 \\
0 & U_{g-\ell'} Z_{g-\ell'}(\tau) {}^t U_{g-\ell'} 
\end{matrix}\right) = g\circ h \left(\begin{matrix}
\tau \cdot 1_\ell & 0 \\
0 & Z_{g-\ell} (\tau) 
\end{matrix}\right).
$$ 
Since $g\circ h$ is an isometry and since $||Z_{g-\ell}(h^{-1}(\tau))|| < |\tau|$, this implies first of all that  
$\ell=\ell'$ and that we can choose $i=i':\Delta^{\ell'} \times \sD_{g-\ell'} \to \sD_g$. Secondly $g\circ h$ and hence $h$ leave the image of $\Delta^\ell$ invariant, as well as the decomposition $\Delta^\ell\times \sD_{g-\ell}\subset \sD_g$. So $h\in G$, as claimed. The permutation of the different 
factors of $i_0(\Delta^\ell)$ defines a homomorphism $G \to \sS_g$ whose kernel is the image 
of $\hat{i}$. We let $H_0$ denote the intersection of $H$ with the image of $\hat{i}$. It is a 
subgroup of $H$ of finite index with the claimed properties. 
\end{proof}
\par
Next we will use that $Y_0$ is a Kobayashi geodesic to construct a splitting
of the map $\varphi_0^*\Omega^1_{\sA_g}\to \Omega^1_{Y_0}$. Equivalently, 
writing $\pi: \widetilde{Y}_0 \to Y_0$ for the universal covering and 
identifying $\widetilde{Y}_0$ with its image in $\sD_g$ we will construct 
a $\pi_1(Y_0,*)$-invariant invertible subsheaf of 
$\Omega^1_{\sD_g}|_{\widetilde{Y}_0}=\pi^*\Omega^1_{\sA_g}$ which splits the 
restriction map $\Omega^1_{\sD_g}|_{\widetilde{Y}_0}\to 
\Omega^1_{\widetilde{Y}_0}$.

\begin{construction}\label{subsheaf}
We keep the notations from Corollary \ref{invariant} and Proposition \ref{interprkob}, c). Choose coordinates 
$\tau_{ij}$ on $\sD_g$ and hence
$\tau_{11},\ldots,\tau_{\ell\ell}$ for the disc $\Delta^\ell$. 
Then $\Omega^1_{\Delta^\ell}=\left\langle d\tau_{11},\ldots,
dx_{\ell\ell} \right\rangle_{\sO_{\Delta^\ell}}$ is a quotient sheaf of $\Omega^1_{\sD_g}|_{\Delta^\ell}$. At the same time, the projections $p_i$ allow to consider $\Omega^1_{\Delta^\ell}$ as a subsheaf of $\Omega^1_{\sD_g}|_{\Delta^\ell}$.

On $\widetilde{Y}_0$ the diagonal subsheaf $\Omega'=\left\langle d\tau_{11}+\cdots+d\tau_{\ell\ell} \right\rangle_{\sO_{\widetilde{Y}_0}}$ of $\Omega^1_{\sD_g}|_{\widetilde{Y}_0}$ satisfies:
\begin{enumerate}
\item[1.] $\Omega'$ is invariant under the action of $H$ and hence $\pi(Y_0,*)$ on
$\Omega^1_{\sD_g}|_{\widetilde{Y}_0}$.
\item[2.] The composition with the natural restriction map 
$$
\Omega^1_{\sD_g}|_{\widetilde{Y}_0} \longrightarrow \Omega^1_{\widetilde{Y}_0}
$$
is an isomorphism, compatible with the action of  $H$ and hence of 
$\pi(Y_0,*)$.  
\end{enumerate}
The first condition implies that $\Omega'$ descends to an invertible subsheaf $\Omega_0$ of
$\Omega_{\sA_g}|_{Y_0}$. Together we obtain:
\begin{enumerate}
\item[3.] $\Omega'=\pi^*\Omega_0$ for an invertible subsheaf
$\Omega_0$ of $\varphi_0^* \Omega_{\sA_g}$ and the composition 
$$
\Omega_0 \subset \varphi_0^* \Omega^1_{\sA_g} \longrightarrow \Omega^1_{Y_0}
$$
is an isomorphism.  
\end{enumerate}
\end{construction}
This completes the proof of the implication ``a) $\Longrightarrow$ b)'' in Theorem~\ref{main} in the case $Y=Y_0$.

\begin{proof}[Proof of ``c) $\Longrightarrow$ a)'' in Theorem~\ref{main}]\ \\
By Lemma~\ref{arakeq} d.) we know that $\V = \BL \otimes \U$. At some point $y \in Y_0$
fix a basis $\{a,b\}$ of the fiber $\BL_y$, a basis $\{e_1,\ldots,e_s\}$
of the fiber $\U$. Write $a_i = a \otimes e_i$, $b_i = b \otimes e_i$,
and choose a symplectic basis $\{a_{s+1}, b_{s+1},\ldots,a_g,b_g\}$ of 
the orthogonal complement $\V^\perp$ in $\W$. We identify the universal cover
of $\sA_g$ in this proof with $\BH_g$. There, the basis extends uniquely to 
a basis of sections of $R^1 \tilde{f}_* \Q_A$, where $\tilde{f}: A 
\to \BH_g$ is the universal family of abelian varieties over $\BH_g$.
\par
As usual in the construction of a period matrix for an abelian variety, there
is a unique basis $\{s_1,\ldots,s_g\}$ of sections of $\tilde{f}_* 
\Omega^1_{A/\BH}$ such that 
$$ \int_{b_i} s_j(\tau) = \delta_{ij} \quad \text{for} \quad \tau \in \BH, 
\quad i,j=1,\ldots,g.$$
Then the matrix $Z$ with entries
$$ Z_{ij} = \int_{a_i} s_j(\tau)  \quad \text{for} \quad \tau\in \BH $$
is the period matrix for the fibre $A_\tau$. The point is that 
${\rm pr}_{11}: Z \mapsto z_{11}$ is the period mapping for $\BL$.
Consequently, ${\rm pr}_{11} \circ \tilde{\varphi}_0$ is an isometry
and by Proposition~\ref{interprkob} the map $\varphi_0$ is a Kobayashi geodesic.
\end{proof}
So we finished the proof of Theorem~\ref{main} for $Y=Y_0$. For the general case
it remains to extend the Construction \ref{subsheaf} to the boundary.
\section{Splitting of the log-differentials}\label{splittings}

We summarize what we need on toroidal compactifications
of $\sA_g$ as used in \cite{mu77}, following the notation of 
loc.~cit.\ whenever possible.
\par
The degree $g'$ of a point $Z$ on the boundary of $\sD_g$ is defined as
the rank of $1-ZZ^*$. Up to the action of $\Sp(2g,\R)$ a point
of degree $g'$  is equivalent to a point in a standard
boundary component
$$ F_{g'} = \left\{\left(\begin{matrix} 1_{g-g'} & 0 \\
0 & Z\\ \end{matrix} \right), \quad Z \in \sD_{g'} \right\} \subset \sD_{g'}.$$
The sets $M \cdot F_{g'}$ for $0 \leq g' <g$ and $M \in \Sp(2g,\R)$
are called {\em boundary components}. A boundary component is called
{\em rational} if $M \in \Sp(2g,\Q)$, provided that the
identification of the universal covering $\widetilde{\sA_g} \cong \sD_g$
has been done using a rational basis of the (relative) homology.
The stabilizer of $F_{g'}$ is the group
\begin{equation} \label{eq:NF}
N(F_{g'}) = \{M \in \Sp(2g,\R), MF_{g'} = F_{g'}\} = 
\left\{ \left(
\begin{matrix} u & * & * & * \\ 
0 & A'  & * & B' \\
0 & 0 & (u^t)^{-1} & 0\\
0 & C' &  * & D' \\
\end{matrix} 
\right)\right\} \subset \Sp(2g,\R)
\end{equation}
We denote the unipotent radical of $N(F)$ by $W(F)$ and
denote the center of $W(F)$ by $U(F)$. See \cite{na80}
for explicit matrix realizations and examples. 
\par
A neighborhood of a boundary component $F$ 
$$D_F = \bigcup_{M \in U(F)_\C} M \cdot F $$
admits Siegel domain coordinates, given for a standard
boundary component $F=F_{g'}$ by
\begin{align*}
 D_{F_{g'}} &\cong U(F_{g'})_\C \times (W(F_{g'})/U(F_{g'}))_\C \times 
F_{g'} \\
\left(\begin{matrix} Z'' & Z''' \\
(Z''')^t & Z'\\ \end{matrix} \right) &\mapsto 
(Z'',Z''',Z') \\
\end{align*}
A toroidal compactification of $\sA_g$ is constructed by gluing $\sA_g$
with toroidal compactifications of the quotients of 
rational boundary components $F$ by $U(F)$. We won't need the details.
\par
Write $x_{ij}$ ($i \leq j$) for the coordinate functions of  
$U(F)_\C \cong M^{(g-g')\times (g-g')}(\C)$, write 
$y_{ij}$ for those of $(W(F)/U(F))_\C \cong M^{g' \times (g-g')}$ and 
write $t_{ij}$ ($i \leq j$) for the coordinate functions of 
$F_{g'} \subset M^{g'\times g'}(\C)$.
By \cite[Proposition~3.4]{mu77} the  pullback of sections
of $\Omega^1_{\bsA_g}(\log S_{\bsA_g})$ to 
$\sD_g$ in the neighborhood a boundary 
component of degree $g'$ to are generated by
$dx_{ij}, dy_{ij}, dt_{ij}$ as $\C[x_{ij},y_{ij},t_{ij}]$-module.
\par
\begin{proof}[Proof of ``a) $\Longrightarrow$ b)'' in Theorem~\ref{main}]\ \\
It remains to show 
that the sheaf $\Omega_0$ in the Construction \ref{subsheaf} extends to 
a subsheaf $\Omega^1\subset \varphi^*\Omega_{\bsA_g}(\log S_{\bsA_g})$, 
isomorphic to $\Omega^1_Y(\log S)$.
\par
Fix a point $s \in S \subset Y$. The universal covering
factors as 
$$ \widetilde{Y_0} \cong \Delta \to \Delta^* \to Y_0,$$
where the first map is $z \mapsto q= e^{2\pi i}$ and $s$ corresponds
to $q=0 \in \Delta^*$. A local generator of $\Omega^1_Y(\log S)$
around $s$ is thus given by $dq/q$ or equivalently by $dz$
near in a rational boundary point, say $1$, of $\Delta$.
\par
We thus need to check that the local generator 
$d\tau_{11}+\ldots d\tau_{\ell\ell}$ of $\Omega'$ is in the pullback
of $\Omega^1_{\bsA_g}(\log S_{\bsA_g})$ to $\sD_g$
near $\varphi(s)$. This is a matter of understanding
a base change.
\par
Any two identifications of $\widetilde{\sA_g} \cong \sD_g$
differ by postcomposition with the action of some 
$M \in \Sp(2g,\R)$. This applies in particular to the 
normalization of the multidisc in Proposition~\ref{interprkob} and
the one where the $F_{g'}$ are rational boundary
components. 
\par
Since $\tilde{\varphi}$ is the
diagonal embedding onto the multidisc of dimension $\ell$
times a map $Z_{g-\ell}(\tau)$ with $||Z_{g-\ell}(\tau)|| < |\tau|$ , 
the point $\varphi(s)$ lies on a boundary component of degree
$g' =  g-\ell$, i.e.\ 
$$ \widetilde{\varphi}(1) \ = \left(\begin{matrix} 1_{\ell} & 0 \\
0 & Z_{g-\ell}(1)\\ \end{matrix} \right) \in \overline{D_g},$$
still with respect to the identification $\widetilde{\sA_g} \cong \sD_g$
using conjugation by the matrix $U$ of Proposition~\ref{interprkob}.
\par
This implies that the base change $M \in N(F_{g'})$. Recall the action of
a symplectic matrix on $D_g$ is given by
$$ \left(\begin{matrix} A & B \\
C & D \\ \end{matrix} \right) \cdot Z = \frac{(A-IC)(Z+1) + (B-ID)I(Z-1)} 
{(A+IC)(Z+1) + (B+ID)I(Z-1)}$$
Using this, one quickly calculates, writing $M \in N(F_{g'})$
in blocks as in equation~\eqref{eq:NF}, 
$$ (dM^*) \left(\begin{matrix} (d\tau_{ij}) & 0 \\
0 & 0 \\ \end{matrix} \right) = 
\left(\begin{matrix} u\cdot (d\tau_{ij}) & 0 \\
0 & 0 \\ \end{matrix} \right).$$
In particular, $(dM^*)(d\tau_{11}+\ldots d\tau_{\ell\ell})$ is
a $\C$-linear combination of $dx_{ij}$, what we needed to show.
\end{proof}
\section{Affine geodesics for the Kobayashi metric}\label{affine}
In this section we will assume that $Y_0$ is affine, hence that $S\neq \emptyset$, and we will show that the image of $Y_0$ in $\sA_g$ lies in a Hilbert modular surface $U$. If the general fibre of $f_0$ is simple, $U$ will be indecomposable, in the sense that no finite \'etale cover can be a product. 
\begin{theorem}\label{HM}
Let $Y_0$ be an affine curve, and let $\varphi_0:Y_0\to \sA_g$ be a geodesic for the Kobayashi metric, induced by a family $f_0:X_0\to Y_0$ of polarized abelian varieties. Assume that the induced variation of Hodge structures $\W_\Q=R^1f_*\Q_{X_0}$ is irreducible, or equivalently that the general fibre of $f_0$ is simple.

Then there exists a totally real number field $\sigma: K\subset \C$ of degree $g$ over $\Q$ and a rank $2$ variation of Hodge structures $\BL_K$ of weight one, such that
\begin{enumerate}
\item $\BL=\BL_K\otimes_K\C$ is maximal Higgs. In particular the Higgs bundle of $\BL$ is of the form 
$$
\begin{CD}
\big(\sL\oplus \sL^{-1}, \tau: \sL@> \cong >> \sL^{-1}\otimes \Omega_Y^1(\log S)\big),    
\end{CD}
$$
where $\sL$ is a logarithmic theta characteristic and $\tau$ the induced isomorphism.
\item If $\sigma_1= \sigma, \sigma_2, \ldots, \sigma_g: K \to \C$ are the different embeddings of $K$ in $\C$, and if $\BL_i=\BL_K^{\sigma_i}$ is the corresponding conjugate of $\BL$, then
$$
\W=\W_\Q\otimes_\Q \C=\BL_1\oplus \cdots \oplus \BL_g.
$$
\item None of the $\BL_i$ is unitary.
\end{enumerate}
\end{theorem}
\begin{corollary}
In Theorem \ref{HM} there exists a $g$-dimensional Hilbert modular 
variety $Z_0\subset \sA_g$
such that $\varphi_0:Y_0\to \sA_g$ factors through $Z_0\subset \sA_g$.
We have $\ell = g$ in the notation of Proposition~\ref{interprkob} 
and the corresponding polydisc is  the universal covering of $Z$.
\end{corollary}
\begin{proof}[Proof of Theorem \ref{HM}]
Several of the arguments used in the proof are taken
from \cite[Section 3 and 4]{VZ04}.
\par
By Theorem \ref{main} and by Lemma \ref{arakeq} $\W$ contains a direct factor $\V$
which is of the form $\BL\otimes \U$ with $\U$ unitary and with
$\BL$ maximal Higgs of rank two, hence induced by a logarithmic theta characteristic $\sL$.

The local system $\U^{\otimes 2}$ is a subsystem of $\W^{\otimes 2}$. So on one hand, the local monodromy in $s\in S$ is unipotent, on the other hand it is unitary, 
of finite order. Then the local residues of $\U^{\otimes 2}$ and hence of $\U$ in $s\in S$ are trivial and $\U$ extends to a unitary local system on $Y$, again denoted by $\U$. Writing
$\sU=\U\otimes \sO_Y$ the residues of $\BL\otimes \U$ in $s\in S$ are given by
$$
\begin{CD}
\sL\otimes \sU @> \tau >> \sL^{-1}\otimes \sU \otimes \Omega_Y^1(\log S)\big)
@> {\rm res} >> \sL^{-1}\otimes \sU|_s,    
\end{CD}
$$ 
hence one finds:
\begin{claim}\label{claim1}
The residues of $\BL\otimes \U$ in $s\in S$ are isomorphisms
$\sL\otimes \sU|_s \to \sL^{-1}\otimes \sU|_s$. 
\end{claim}
Since the rank two variations of Hodge structures, given by a theta characteristics, 
are unique up to the tensor product with local systems induced by two-division points,
we can write
$$
\W = \W_1 \oplus \T \mbox{ \ \ \ with \ \ \ } \W_1=\BL\otimes \U \oplus \V'
$$
where $\T$ is the maximal unitary local subsystem, and where none of the direct
factors of $\V'$ is maximal Higgs. By abuse of notations we write $\V=\BL\otimes \U$,
allowing $\V$ to be reducible.
\begin{claim}\label{claim2}
The unitary part $\T$ is zero. 
\end{claim}
\begin{proof}
Assume the contrary. By \cite[Lemma 3.3]{VZ04} the first decomposition
is defined over a number field. As above the residues of $\T$ in $s\in S$ are zero. This being invariant under Galois conjugation, the same holds true for all conjugates
of $\T$ in $\W$. Since we assumed that $\W$ is irreducible over $\Q$, the same holds 
for $\W$, contradicting Claim \ref{claim1}.
\end{proof}
\begin{claim}\label{claim3}
The direct sum decomposition $\W=\V \oplus \V'$ is defined over a number field. 
\end{claim}
\begin{proof}
If not, as in the proof of \cite[Lemma 3.3]{VZ04}, one obtains a non-trivial family
$\V'_t$ of local systems over the disk $\Delta$ with $\V'_0=\V'$.
The non-triviality of the family forces the composition
$$
\V'_t \to \W \longrightarrow \BL\otimes \U
$$
to be non-zero, for some $t$ arbitrarily close to zero. The complete reducibility of local systems, coming from variations of Hodge structures, implies that some local system
$\V''$ lies in both, $\V'_t$ and $\BL\otimes \U$. Since for the second one the Higgs fields are isomorphism, the same holds true for $\V''$. On the other hand, for $t$ sufficiently small, $\V'_t$ will not contain any local subsystem with isomorphisms as Higgs fields.
\end{proof}
Since $\V$ is defined over a number field, we can apply \cite[Corollary 3.7]{VZ04}
and we find that $\BL$ and $\U$ can be defined over a number field, as well.

Consider next the sublocal system $\V^{\otimes 2}$ of $\W^{\otimes 2}$. Since $\BL$
is self dual, one can identify $\BL^{\otimes 2}$ with $\BEnd (\BL)$, hence one obtains a
trivial direct factor and its complement, denoted by $\BL^{\otimes 2}_0$ in the sequel.
So $\V^{\otimes 2}$ decomposes as a direct sum of $\BL^{\otimes 2}_0 \otimes \U^{\otimes 2}$ and
$\U^{\otimes 2}$. For the first one the residues are again given by isomorphisms,whereas for the second one the residues are zero. 

This property is preserved under Galois conjugation. Since for $\sigma \in {\rm Gal}(\bar{\Q}/\Q)$ one has $\V^\sigma=\BL^{\sigma}\otimes \U^{\sigma}$, and since the weight is zero, either $\BL^{\sigma}$ or $\U^{\sigma}$ must be unitary. The first case can not occur, since the residues are isomorphisms. So the maximal unitary subbundle $\T'$ of
$\W^{\otimes 2}$ is the Weil restriction of $\U^{\otimes 2}$. So it is defined over $\Q$ and inherits a $\Z$ structure from the one of $\W^{\otimes 2}$. This implies that 
$\T'$ trivializes after replacing $Y_0$ by an \'etale covering, and hence $\U^{\otimes 2}$ as well. In particular writing $\U=\U_1\oplus \cdots \oplus \U_s$ for the decomposition of $U$ in $\C$-irreducible direct factors, one finds
$$
\U^{\otimes 2}=\bigoplus_{i,j} \U_i\otimes \U_j\cong \bigoplus \C,
$$ 
which is only possible if all $\U_j$ are isomorphic and of rank one, hence trivial.\vspace{.2cm}

Up to now we verified the existence of a direct sum decomposition of $\W$, satisfying the conditions (1), (2) and (3), over some number field $L$. For the next step, we will not need that $Y_0$ is affine. 
\begin{claim}\label{claim5}
Assume that over some curve $Y_0$ and for some number field $L$, there exists a direct sum decomposition, satisfying the conditions (1), (2) and (3) stated in Theorem
\ref{HM}, and assume that $L$ is minimal with this property. Then $[L:\Q]=g$, and $L$ is either
totally real, or an imaginary quadratic extension of a totally real number field $K$. 
\end{claim}
\begin{proof}
Let $L$ be the field of moduli of $\BL_1$, i.e. $L$ is the
fixed field of all $\sigma \in \Gal_{\overline{\Q}/Q}$
such that $\BL_1^\sigma \cong \BL_1$. Obviously one finds that $[L:\Q]$ is equal to the number of conjugates, hence equal to $g$.

Let $K$ be the trace
field of $\BL_1$. Since $\BL_1$ is isomorphic to a Fuchsian
representation, $K$ is real. It is well known (e.g. \cite{Ta69}) that a rank two
local system can be defined over an extension of degree at most two over its trace field. In particular,
$[L:K] \leq 2$.
\par
We claim that $L$ injects into the endomorphism ring of
the family $f_0$. In fact, let $F$ be the Galois closure of $L/\Q$ and
choose coset representatives $\tau_i$ of $\Gal(F/\Q)/\Gal(F/L)$.
Then for $x \in F$,
$$ \sum_{\sigma \in \Gal(F/\Q)} \sigma(x)\cdot id_{\BL^\sigma}
= \sum_{\tau_i} \tau_i({\rm tr}^F_L(x)) \cdot id_{\BL^{\tau_i}}$$
is in $\End_\Q(\W_Q)$ of bidegree $(0,0)$, hence an
endomorphism of $f_0$. Specializing to $x \in L$ proves the claim.
\par
By Albert's classification (\cite{BL04} \S 5.5) of endomorphism rings
and since the general fibre of $f_0$ is simple,
$L$ is totally real or an imaginary quadratic extension
of a totally real field $K$.
\end{proof}
To finish the proof of Theorem \ref{HM} it remains to show that 
the minimal field $L$ in Claim \ref{claim5} has to be totally real, if
$Y_0$ is affine. If $K=L$, there is nothing left to show.
\par
We distinguish cases according to what the full endomorphism ring
$W_\Q$ is, following notations  in \cite{BL04} \S 5.5. 
First, suppose the endomophism ring is of the second kind, i.e.\ the Rosati involution fixes $L$. 
Then since $[K:\Q] = \rk(\W)$ and since all the $\BL_i$ have
non-trivial Higgs field, the endomorphism ring is in fact an 
indefinite quaternion algebra $B$, containing $L$ (\cite{BL04} 9.10.Ex.~(4), 
\cite{Sh63} Proposition~18).
\par
Second, suppose that the endomorphism ring is a totally definite quaternion
algebra. Then (by \cite{BL04} 9.10.Ex~(1) or \cite{Sh63} Proposition~15)
the general fibre of $f_0$ is not simple.
\par
It remains thus only to discuss the third case, a totally indefinite quaternion
algebra $B$. Let $Z$ be the (``PEL''-) Shimura variety of abelian varieties
with $B\subset \End(\W_\Q)$. To construct $Z$, let $G$ be the $\Q$-algebraic
group with $G(\Q) = B^*$, fix a maximal compact subgroup $K$ and form the quotient of $G/K$
by an arithmetic lattice $\Gamma$. For the universal family $A$ over $\Gamma\backslash G/K$ we
have, by construction, a homomorphism $B \to \End_Q(A)$, injective since $B$ is simple.
For dimension reasons we conclude that $Z = \Gamma\backslash G/K$ and that $Z$ is compact
by \cite{Sh71} Proposition~9.3. Thus the map $Y_0\to Z$ would extend to $Y\to Z$, 
contradicting the maximality of the Higgs field.
\end{proof}
\begin{remark} If in Claim \ref{claim5} $L\neq K$, then $Y_0$ is projective. For $K=\Q$ examples are the moduli schemes of {\em false elliptic curves,} i.e. of polarized abelian surfaces $Z$ with ${\rm End}(Z)\otimes_\Z\Q$ a totally indefinite quaternion algebra over $\Q$. 

In a similar way there should exist compact {\em false Teichm\"uller curves} in $\sA_g$.
The corresponding variation of Hodge structures $\W_\Q$ should have a $K$ irreducible direct factor $\T_K$ of rank $4$, for some totally real number field $K$, and
$\T=\T_K\otimes \C$ should be the direct sum of two complex variations of Hodge structures of rank $2$.
\end{remark}
\section{Families of Kobayashi geodesics}\label{fam}
Let $k$ be an algebraically closed subfield of $\C$, let $U$ be an irreducible
non-singular variety, defined over $k$ and let $h:\sY \to U$ be a smooth family of projective curves. We consider a relative normal crossing divisor $\sS\in \sY$ and write $\sY_0=\sY\setminus \sS$.
 
We assume in the sequel that $\tilde{\varphi}_0:\sY_0 \to \sA_g$ is a generically finite morphism, defined over $k$, which extends to a morphism $\tilde{\varphi}:\sY\to \bsA_g$. Here we consider $\bsA_g$ as a $k$-variety. So $\tilde{\varphi}$ and $h$ define a morphism
$$
\varphi':\sY\to \bsA_{g \, U}:= \bsA_g\times U.
$$ 
We want that all fibres of $h|_{\sY_0}$ are Kobayashi geodesics. However, since we do not want to study the behavior of this property in families, we use instead an extension of condition b) in Theorem \ref{main} to the relative case, hence the natural morphism
\begin{equation}\label{relsplit}
\Psi:\varphi'^*\Omega^1_{\bsA_{g \, U}/U}(\log(S_{\bsA_{g \, U}})) \longrightarrow
\Omega^1_{\sY/U}(\log \sS ),
\end{equation}
where $S_{\bsA_{g \, U}}=S_{\bsA_g}\times U = \bsA_{g \, U}\setminus \sA_g\times U$.
\begin{lemma}\label{families}
We keep the notations introduced above and assume that the morphism $\Psi$ in
\eqref{relsplit} splits, that $\tilde{\varphi}:\sY \to \bsA_g$ is generically finite
and that the general fibre of the pullback to $\sY$ of the universal family of abelian varieties is simple. Then
\begin{enumerate}
\item[i.] The family $\sY_0 \to U$ is locally trivial, i.e. locally in the analytic topology a product. 
\item[ii.] If for $u\in U$ in general position the fibre $Y=h^{-1}(u)$ is affine, then $U$ is a point.
\end{enumerate}
\end{lemma}
\begin{proof}
The splitting of the map $\Psi$ in \eqref{relsplit} implies that on all fibres $Y$ of $h$
one has a splitting of
\begin{equation}\label{split}
\varphi^*\Omega^1_{\bsA_g}(\log S_{\bsA_g}) \longrightarrow
\Omega^1_{Y}(\log S ).
\end{equation}
So by Theorem \ref{main} all fibres of $h$ are Kobayashi geodesics.
Let $\Omega$ denote a subsheaf of $\varphi'^*\Omega^1_{\bsA_{g \, U}/U}(\log(S_{\bsA_{g \, U}}))$ for which 
$$
\Psi|_{\Omega}:\Omega \to\Omega^1_{\sY/U}(\log \sS )
$$
an isomorphism. Since  
$$
\varphi'^*\Omega^1_{\bsA_{g \, U}}(\log(S_{\bsA_{g \, U}}))=\varphi'^*\Omega^1_{\bsA_{g \, U}/U}(\log(S_{\bsA_{g \, U}}))\oplus h^*\Omega^1_U
$$
the subsheaf $\Omega\oplus \{0\}$ defines a splitting of the composition
$$
\varphi'^*\Omega^1_{\bsA_{g \, U}}(\log(S_{\bsA_{g \, U}}))\longrightarrow
\varphi'^*\Omega^1_{\bsA_{g \, U}/U}(\log(S_{\bsA_{g \, U}})) \longrightarrow
\Omega^1_{\sY/U}(\log \sS ).
$$
This composition factors through
$$
\varphi'^*\Omega^1_{\bsA_{g \, U}}(\log(S_{\bsA_{g \, U}}))\longrightarrow
\Omega^1_\sY(\log S_{\sY}) \longrightarrow
\Omega^1_{\sY/U}(\log \sS ).
$$
Hence the exact sequence
\begin{equation}\label{exsesp}
0\to h^* \Omega^1_U \to \Omega^1_\sY(\log S_{\sY}) \to
\Omega^1_{\sY/U}(\log \sS ) \to 0
\end{equation}
splits as well. Recall that the Kodaira-Spencer map
$$
\rho_u:T_{u,U} \longrightarrow H^1(Y, T_Y(-\log S)),
$$
controlling the infinitesimal deformations of $Y\setminus S$ for $S=\sS_u$ and $Y=h^{-1}(u)$, is given by edge morphism of the dual exact sequence of \eqref{exsesp},
restricted to $u\in U$ (see \cite{Ka78}, for example). So the splitting of
\eqref{exsesp} implies that $\rho_u$ is zero, for all $u\in U$.
By \cite[Corollary 4]{Ka78} this implies that the family $\sY_0\to U$
is locally a product, as claimed in i).

From now on we will use the analytic topology, hence assume by abuse of notations that
$\sY=Y_0\times U$. Part ii) of the Lemma will follow, if we show that
the morphism $\varphi_0:Y_0 \to \sA_g$ is rigid. Using the description of the infinitesimal deformations of subvarieties of $\sA_g$ in \cite{Fa83},
one has to verify that $\End(R^1f_*\Q_{X})=\End(\W_\Q)$ is concentrated in bidegree $(0,0)$. 

Remark under the assumption in ii) Theorem \ref{HM} allows to write
$$
\W=\BL_1^{\oplus \nu_1} \oplus \cdots \oplus \BL_\ell^{\oplus \nu_\ell},
$$
where $\BL_1=\BL$ satisfies the Arakelov equality,
where $\BL_1, \ldots , \BL_\ell$ are pairwise non-isomorphic, and where
all the $\BL_j$ are irreducible of rank $2$. 
Then there are no global homomorphisms $\BL_i \to \BL_j$, except for $i=j$, and 
$\End(\BL_j)^{0,0}=\End(\BL_j)\cong \C$. Hence  
$$
\End(\W)= \bigoplus_{j=i}^\ell M(\nu_j,\End(\BL_j))= M(\nu_j,\C). 
$$
is equal to $\End(\W)^{0,0}$.
\end{proof}
\begin{corollary}\label{field}
Let $\varphi_0:Y_0\to \sA_g$ be an affine Kobayashi geodesic, such that the
induced variation of Hodge structures $\W_\Q$ is $\Q$-irreducible.
Then $\varphi_0:Y_0 \to \sA_g$ can be defined over a number field.   
\end{corollary}
\begin{proof}
Let $\bsA$ be a Mumford compactification, defined over $\bar{\Q}$. Recall that we assumed
$\sA_g$ to be a fine moduli scheme, and that $Y_0$ is a Kobayashi geodesic, if its image in $\sA_g$ has this property. So we may assume that $\varphi_0:Y_0 \to \sA_g$ is injective. 

The $\C$-morphism $\varphi_0:Y_0 \to \sA_g$ extends to a morphism $\varphi:Y \to \bsA_g$, and by Theorem \ref{main} there exists a subsheaf $\Omega \subset \varphi^*\Omega^1_{\bsA_g}(\log S_{\bsA_g})$ splitting the natural projection 
$$
\varphi^*\Omega^1_{\bsA_g}(\log S_{\bsA_g})\to \Omega^1_Y(\log S).
$$
Choose a field $K\subset \C$, finitely generated over $\bar{\Q}$ such that
$Y$, $S$, $\varphi:Y\to \bsA_g$, and $\Omega$ are defined over $K$. By abuse of notations we will use the same letters for the objects, as schemes, morphisms or sheaves over $K$.

Let $U$ denote a non-singular quasi-projective $\bar{\Q}$-variety with $K=\bar{\Q}(U)$.
Choosing $U$ small enough, we can assume (step by step):
\begin{itemize}
\item $Y\to {\rm Spec}(K)$ is the general fibre of a projective morphism
$h:\sY \to U$.
\item $\varphi: Y\to \bsA_g$ extends to a morphism $\tilde{\varphi}:\sY \to \bsA_g$.
\item $h$ is smooth and $\sS:=\tilde{\varphi}^{-1}(S_{\bsA_g})_{\rm red}$ consists of
disjoint sections of $h$.
\item $\Omega \subset \varphi^*\Omega^1_{\bsA_g}(\log S_{\bsA_g})$ is induced by a subsheaf $\tilde{\Omega }\subset \tilde{\varphi}^*\Omega^1_{\bsA_g}(\log S_{\bsA_g})$.
\item $\tilde{\Omega }$ defines a splitting of the natural map
$$
\tilde{\varphi}^*\Omega^1_{\bsA_g}(\log S_{\bsA_g})=
(\tilde{\varphi}\times h)^*\Omega^1_{\bsA_g\times U/U}(\log (S_{\bsA_g}\times U))
\longrightarrow \Omega^1_{\sY/U}(\log \sS). 
$$
\end{itemize}
Writing $\sY_0=\sY\setminus \sS$ we claim that $\dim(\tilde{\varphi}(\sY_0))=1$. 
Otherwise we may replace $U$ by the intersection of $\dim(\tilde{\varphi}(\sY_0))-2$ general hyperplane section, hence assume that besides of the properties stated above 
$U$ is a curve and the morphism $\tilde{\varphi}$ generically finite, 
contradicting the Lemma \ref{families}.

If $\dim(\tilde{\varphi}(\sY_0))=1$ choose any $\bar{\Q}$-valued point in $U$.
Then 
$$
\tilde{\varphi}(h^{-1}(u)\cap \sY_0)= \tilde{\varphi}(\sY_0)
$$ 
and $Y_0$ is obtained by base extension from $h^{-1}(u)\cap \sY_0$.
\end{proof}
\begin{remark}\label{remfield}
As stated in Lemma \ref{families} the proof of Corollary \ref{field} is based on two observations. Firstly the family $\sY_0\to U$ constructed in the proof of the Corollary is locally a product, and secondly one has rigidity of the embedding of a fixed curve as a Kobayashi geodesics. 

In \cite{Moe04} Corollary \ref{field} was formulated for Teichm\"uller curves, but the first part of the argument was omitted.  
\end{remark}
\begin{remark}\label{remunitary}
Without a bit of classification of the possible variations of Hodge structures we are not able to give a criterion for rigidity of compact Kobayashi geodesics. As for Shimura curves the existence of a large unitary subbundle in $\W$ implies non-rigidity.
However there might be other criteria. 
\end{remark}

\section{Examples and comments}\label{examples}

Let $K = \Q({\rm tr}(\gamma), \gamma \in \Gamma)$ be the trace field of
the uniformizing group $Y_0 = \BH /\Gamma$. Moreover 
let $d = [K:\Q]$ and $B = K\cdot \Gamma$ be the subalgebra of ${\rm GL}_2(\C)$
generated by $\Gamma$.
\par
Among the curves, satisfying the equivalent conditions in Theorem~\ref{main}, one finds:\\[.2cm] 
{\bf i)} Shimura curves (\cite{VZ04}). In this case the $\C$ local system $\W$ decomposes as a direct sum $\V_1\oplus \cdots \oplus \V_\ell$ and each $\V_i$ is either unitary or maximal Higgs. 
The algebra $B$ is unramified at all but one of the infinite places of $K$. The curve $Y_0$ can be both compact and non-compact in this case.\\[.1cm]
{\bf ii)} Teichm\"uller curves. By \cite{VZ06} and \cite{Moe04} 
these are the  Kobayashi geodesic curves $Y_0\subset \sA_g$ that 
lie in the image of $M_g$. In this case $\ell=1$ by \cite{Moe04}
and if $d>1$ then $\V^\perp$ is not unitary. $B$ is ramified at
all the infinite places of $K$. The curve $Y_0$ is not compact in this case.\\[.1cm]
{\bf iii)}  There are examples of compact Kobayashi geodesics $Y_0$
that are not Shimura curves. In fact, by \cite{DeMo86} (see also 
\cite{CoWo90})  all triangle groups $\Gamma=\Delta(l,m,n)$ arise as
uniformizing groups of  Kobayashi geodesic curves $Y_0\subset \sA_g$.
If $l,m,n < \infty$ are chosen  not to be in Takeuchis's finite list of arithmetic
triangle groups (\cite{Ta75}), then $Y_0$ is neither a Shimura curve nor a Teichm\"uller curve.\\[.1cm]
{\bf iv)} Not all Kobayashi geodesics in $\sA_g$ with $d>1$, 
with $Y_0$ not compact and $B$ ramified at all the infinite places 
of $K$ are Teichm\"uller curves in $M_g$. In fact, the Prym construction
of \cite{McM06} gives Teichm\"uller curves in $M_4$, whose family
of Jacobian splits. The Prym-antiinvariant part together with $1/2$
the pullback polarization gives such a Kobayashi geodesic in $A_2$.\\[.2cm]
Given this list of examples, the following questions remain open:
\begin{itemize}
\item[1)] Are there compact Kobayashi geodesics in $\sA_g$, 
that are neither Shimura curves nor uniformized by a group
commensurable to a triangle group?
\item[2)] Can one classify all non-Shimura Kobayshi geodesics on the image of a Hilbert 
modular surface in $A_2$ for a given discriminant?
\end{itemize}

\end{document}